\documentclass{llncs}

\usepackage{graphicx}

\usepackage{tikz}

\usepackage{color}

\usepackage{subfigure} 

\usepackage{amsmath}   

\usepackage{amssymb}

\usepackage{float}

\usepackage{ifpdf}

\newcommand{\Cauchy}{L}

\newcommand{\IR}{\mathbb{R}}

\begin{document}

\frontmatter          

\title{ PDE-constrained LDDMM via geodesic shooting and inexact Gauss-Newton-Krylov optimization using the incremental adjoint Jacobi equations }

\titlerunning{PDE-constrained LDDMM ...}     
%
\author{Monica Hernandez}

\institute{Computer Sciences Department \\ Aragon Institute on Engineering Research \\ University of Zaragoza \\ mhg@unizar.es}

\maketitle

\begin{abstract}

The class of non-rigid registration methods proposed in the framework of PDE-constrained 
Large Deformation Diffeomorphic Metric Mapping is a particularly interesting family of physically 
meaningful diffeomorphic registration methods. 
Inexact Newton-Krylov optimization has shown an excellent numerical accuracy and an extraordinarily 
fast convergence rate in this framework.
However, the Galerkin representation of the non-stationary velocity fields does not provide proper geodesic paths.
In this work, we propose a method for PDE-constrained LDDMM parameterized in the space of initial velocity fields
under the EPDiff equation.
The derivation of the gradient and the Hessian-vector products are performed on the final velocity field and
transported backward using the adjoint and the incremental adjoint Jacobi equations.
This way, we avoid the complex dependence on the initial velocity field in the derivations and the 
computation of the adjoint equation and its incremental counterpart.
The proposed method provides geodesics in the framework of PDE-constrained LDDMM, and 
it shows performance competitive to benchmark PDE-constrained LDDMM and EPDiff-LDDMM methods.

\end{abstract}


\section{Introduction}

In the last two decades, diffeomorphic image registration has arisen as a powerful paradigm for deformable image registration~\cite{Sotiras_13}. 
Diffeomorphic registration methods compute transformations fundamental in Computational Anatomy applications~\cite{Kobatake_17}.
Although the differentiability and invertibility of a diffeomorphism constitute fundamental properties for its use in Computational Anatomy,  
the diffeomorphic constraint does not necessarily guarantee that a transformation computed with a given method is physically meaningful for the clinical domain of interest. 
PDE-constrained diffeomorphic registration methods constitute an appealing paradigm for computing transformations under 
plausible physical models of interest.

The first class of diffeomorphic image registration methods constrained to PDEs arose with the work of Younes et al.~\cite{Younes_07},
followed by the proposals in~\cite{Ashburner_11,Vialard_11,Zhang_15}. 
In these works, the transformations are parameterized by time-varying velocity fields that satisfy the EPDiff equation~\cite{Holm_98}. 
This physical constraint allows formulating the problem in the space of initial velocity fields. 
This guarantees that the obtained transformations belong to geodesic paths of diffeomorphisms,
which is desirable in important Computational Anatomy applications~\cite{Miller_05}.
For this class of methods, the dependence of the energy functional on the initial velocity field is complex, 
and the optimization is usually implemented using gradient-descent. 

The class of diffeomorphic registration methods proposed in~\cite{Hart_09,Vialard_11,Mang_15,Mang_16,Mang_17} is especially interesting,
where the physical PDEs are imposed using hard constraints. 
These methods model the problem using a PDE-constrained variational formulation.
Numerical optimization is approached using gradient-descent~\cite{Hart_09,Vialard_11,Mang_15} 
and second-order optimization in the form of inexact reduced Newton-Krylov methods~\cite{Mang_15,Mang_16,Mang_17}.  
The gradient and the Hessian-vector products are computed by the differentiation of the augmented energy functional using optimal control theory methods.
The constrained optimization approach provides the versatility to impose different physical models to the computed transformations 
by simply adding the PDEs associated to the problem as hard constraints.
In particular, inexact Newton-Krylov optimization shows an excellent numerical accuracy and an extraordinarily fast convergence rate~\cite{Mang_15,Mang_16,Mang_17}.
However, 
the Galerkin representation of the non-stationary velocity fields used for diffeomorphism parameterization 
does not provide proper geodesic paths.

The purpose of this article is to parameterize the compressible method for PDE-constrained diffeomorphic registration with inexact Newton-Krylov optimization in~\cite{Mang_16} 
in the space of initial velocity fields.
Rather than computing the gradient and the Hessian-vector products from the differentiation of the augmented energy functional on the initial velocity field,
we perform the computations on the final velocity field and transport the computations backward using the adjoint and the incremental adjoint Jacobi equations.
This way, we avoid the complex dependence on the initial velocity field in the derivations.
The proposed method provides geodesics in the framework of PDE-constrained LDDMM.
The method is evaluated and compared to benchmark PDE-constrained LDDMM~\cite{Mang_15} and EPDiff-LDDMM~\cite{Zhang_15} methods with the 
manual segmentations of the Non-Rigid Image Registration Evaluation Project (NIREP) database~\cite{Song_10}.


In the following, Section~\ref{sec:BackgroundMethods} revisits the methods more related to our approach. 
Section~\ref{sec:Method} presents our proposed method.
Section~\ref{sec:Results} shows the most relevant experiments for the evaluation of our method.
Finally, Section~\ref{sec:Conclusions} gathers the most remarkable conclusions of our work.


\section{Background methods}
\label{sec:BackgroundMethods}

\subsection{Large Deformation Diffeomorphic Metric Mapping (LDDMM)}

Let $I_0$, and $I_1$ be the source and the target images defined on the image domain $\Omega \subseteq \mathbb{R}^d$.
We denote with $Diff(\Omega)$ to the Riemannian manifold of diffeomorphisms on $\Omega$.
$V$ is the tangent space of the Riemannian structure at the identity diffeomorphism, $id$.
$V$ is made of smooth vector fields on $\Omega$.
The Riemannian metric is defined from the scalar product in $V$ 
\begin{equation}
 \langle v, w \rangle_V = \langle \Cauchy v, w \rangle_{L^2} = \int_\Omega \langle \Cauchy v(x), w(x) \rangle d\Omega, 
\end{equation}
\noindent where $\Cauchy = (Id - \alpha \Delta)^s, \alpha >0, s \in \mathbb{N}$
is the invertible self-adjoint differential operator associated with the differential structure of $Diff(\Omega)$.
We denote with $K$ to the inverse of $L$.

The LDDMM variational problem is given by the minimization of the energy functional
\begin{equation}
\label{eq:LDDMM}
E(v) =  \int_0^1 \langle \Cauchy v_t, v_t \rangle_{L^2} dt + \frac{1}{\sigma^2} \Vert I_0 \circ (\phi^{v}_{1})^{-1} - I_1\Vert_{L^2}^2. 
\end{equation}

\noindent The problem is posed in the space of time-varying smooth flows of velocity fields in $V$, ${v} \in L^2([0,1],V)$.
Given the smooth flow ${v}:[0,1] \rightarrow V$, $v_t:\Omega \rightarrow \IR^{d} \in V$, 
the diffeomorphism $\phi_1^{v}$ is defined as the solution at time $1$ to the transport equation 
$d_t \phi_t^{v} = v_t \circ \phi_t^{v}$ with initial condition $\phi_0^{v} = id$.
The transformation $(\phi^{v}_{1})^{-1}$ computed from the minimum of $E({v})$ is the diffeomorphism that 
solves the LDDMM registration problem between $I_0$ and $I_1$.
The optimization of Equation~\ref{eq:LDDMM} was originally approached using gradient-descent in $L^2([0,1],V)$, yielding the update equation $v_t^{n+1} = v_t^n - \epsilon (\nabla_v E(v))_t$~\cite{Beg_05}.

\subsection{LDDMM in the space of initial velocity fields (EPDiff-LDDMM)}

The geodesics of $Diff(\Omega)$ under the right-invariant Riemannian metric are uniquely determined by 
the time-varying flows of velocity fields that satisfy the Euler-Poincar\'e equation (EPDiff)~\cite{Holm_98}
\begin{equation}
\label{eq:EPDiff}
\partial_t v = -ad_v^\dagger v = - K ad_v^* L v = - K [ (Dv)^T Lv + D(Lv) v + L v \nabla \cdot v ].
\end{equation}
\noindent with initial condition $v_0 \in V$.
LDDMM can be posed in the space of initial velocity fields

\begin{equation}
\label{eq:EPDiffLDDMM}
E(v_0) = \langle \Cauchy v_0, v_0 \rangle_{L^2} + \frac{1}{\sigma^2} \Vert I_0 \circ (\phi^{v}_{1})^{-1} - I_1 \Vert_{L^2}^2, 
\end{equation}

\noindent where $(\phi^{v}_{1})^{-1}$ is the solution at time $1$ to the transport equation of the flow $v_t$ that satisfies 
the EPDiff equation for $v_0$.

The optimization of Equation~\ref{eq:EPDiffLDDMM} was originally approached using gradient-descent in $V$~\cite{Younes_07}
\begin{equation}
 v_0^{n+1} = v_0^n - \epsilon \nabla_{v_0} E(v_0). 
\end{equation}
\noindent More recently, it has been proposed to compute the gradient at $t = 1$ and to integrate backward
the reduced adjoint Jacobi field equations~\cite{Bullo_95} 
\begin{eqnarray}
\label{eq:AdjointJacobiField}
 \partial_t U + ad_v^\dagger U = 0 \textnormal{ in } \Omega \times [0,1) \\
 \partial_t \delta v + U - ad_v \delta v + ad_{\delta v}^\dagger v = 0 \textnormal{ in } \Omega \times [0,1)
\end{eqnarray}
\noindent with initial conditions $U(1) = \nabla_{v_1} E(v_0)$ and $\delta v(1) = 0$,
to get the gradient update at $t = 0$, $v_0^{n+1} = v_0^n - \epsilon \delta v(0)$~\cite{Zhang_15}. 

\subsection{PDE-constrained LDDMM (PDE-LDDMM)}

The variational problem is given by a slight modification of Equation~\ref{eq:LDDMM}
\begin{equation}
\label{eq:LS-LDDMM}
E(v) = \int_0^1 \langle \Cauchy v_t, v_t \rangle_{L^2} dt + \frac{1}{\sigma^2} \Vert m(1) - I_1\Vert_{L^2}^2, 
\end{equation}

\noindent where $m(t) = I_0 \circ (\phi^{v}_{t})^{-1}$ is computed from the solution of the state equation
\begin{equation}
\label{eq:MassTransport}
\partial_t m(t) + \nabla m(t) \cdot v_t = 0 \textnormal{ in } \Omega \times (0,1],
\end{equation}
\noindent with initial condition $m(0) = I_0 \textnormal{ in } \Omega$. 
The combination of Equations~\ref{eq:LS-LDDMM} and~\ref{eq:MassTransport} leads to a PDE-constrained optimization problem.

Optimization can be performed combining the method of Lagrange multipliers with gradient-descent~\cite{Hart_09} or second-order inexact Newton-Krylov methods~\cite{Mang_15}. 
The gradient is computed from  
\begin{eqnarray}
\partial_t m(t) + \nabla m(t) \cdot v_t = 0 \textnormal{ in } \Omega \times (0,1]  \label{eq:FirstOrder1} \\
-\partial_t \lambda(t) - \nabla \cdot ( \lambda(t) \cdot v_t ) = 0 \textnormal{ in } \Omega \times [0,1)  \label{eq:FirstOrder2}  \\
(\nabla_v E(v))_t = \Cauchy v_t + \lambda(t) \cdot \nabla m(t) \textnormal{ in } \Omega \times [0,1] \label{eq:FirstOrder3}  
\end{eqnarray}
\noindent subject to the initial and final conditions $m(0) = I_0$ and $\lambda(1) = -\frac{2}{\sigma^2}(m(1) - I_1)$. 
Equations~\ref{eq:FirstOrder1} and~\ref{eq:FirstOrder2} correspond with the state and adjoint equations, respectively.

The Hessian is computed from
\begin{eqnarray}
\label{eq:SecondOrderVariation}
\partial_t \delta m(t) + \nabla \delta m(t) \cdot v_t + \nabla m(t) \cdot \delta v(t) = 0 \textnormal{ in } \Omega \times (0,1]  \label{eq:SecondOrder1}  \\
-\partial_t \delta \lambda(t) - \nabla \cdot ( \delta \lambda(t) \cdot v_t ) - \nabla \cdot ( \lambda(t) \cdot \delta v(t) ) = 0 \textnormal{ in } \Omega \times [0,1) \label{eq:SecondOrder2}  \\
(H_v E(\delta v))_t = \Cauchy \delta v(t) + \delta \lambda(t) \cdot \nabla m(t) + \lambda(t) \cdot \nabla \delta m(t) \textnormal{ in } \Omega \times [0,1] \label{eq:SecondOrder3} & 
\end{eqnarray}

\noindent subject to the initial and final conditions $\delta m(0) = 0$ and $\delta \lambda(1) = -\frac{2}{\sigma^2}\delta m(1)$.
Equations~\ref{eq:SecondOrder1} and~\ref{eq:SecondOrder2} correspond with the incremental state and incremental adjoint equations, respectively.

The minimization using inexact Newton-Krylov optimization yields to the update equation
\begin{equation}
\label{eq:Newton}
 v_t^{n+1} = v_t^n - \epsilon \delta v_t^n,
\end{equation}
\noindent where $ \delta v^n$ is computed from preconditioned conjugate gradient (PCG) on the reduced system 
\begin{equation}
\label{eq:KKT}
H_v E( \delta v^n) = - \nabla_v E(v^n),
\end{equation}
\noindent with preconditioner $K$.
Since the update equation is written on the reduced space $V$, the minimization is a reduced space optimization method.

By construction, the Hessian is positive definite in the proximity of a local minimum.
However, it can be indefinite or singular far away from the solution.
In this case, the search directions obtained with PCG are not guaranteed to be descent directions. 
In order to overcome this problem, one can use a Gauss-Newton approximation dropping expressions of $H_v E( \delta v^n)$
to guarantee that the matrix is definite positive.
In particular, one can drop the terms in Equations~\ref{eq:SecondOrder2} and~\ref{eq:SecondOrder3} involving the adjoint variable $\lambda$.


\section{The proposed method (PDE-EPDiff LDDMM)}
\label{sec:Method}

The PDE-constrained problem is given by the minimization of the energy functional
\begin{equation}
 E(v_0) = \langle L v_0, v_0 \rangle_{L^2} + \frac{1}{\sigma^2} \Vert m(1) - I_1 \Vert_{L^2}^2,
\end{equation}
\noindent subject to the EPDiff and the state equations
\begin{eqnarray}
\partial_t v_t + ad_{v_t}^\dagger v_t = 0 \textnormal{ in } \Omega \times (0,1] \\
\partial_t m(t) + \nabla m(t) \cdot v_t = 0 \textnormal{ in } \Omega \times (0,1],
\end{eqnarray}

\noindent with initial conditions $v(0) = v_0$ and $m(0) = I_0$, respectively.

Optimization is performed combining the method of Lagrange multipliers with inexact Gauss-Newton-Krylov methods
in the following way. Let
$w: \Omega \times [0,1] \rightarrow \mathbb{R}^d$ and
$\lambda: \Omega \times [0,1] \rightarrow \mathbb{R}$ 
be the Lagrange multipliers associated with the EPDiff and the state equations.
We build the augmented Lagrangian
\begin{multline}
\label{eq:MCC-LDDMM-Aug}
E_{aug}(v_0) =  E(v_0) + \int_0^1 \langle w(t), \partial_t v(t) + ad_{v_t}^\dagger v_t \rangle_{L^2} dt \\
+ \int_0^1 \langle \lambda(t), \partial_t m(t) + \nabla m(t) \cdot v_t \rangle_{L^2} dt.
\end{multline}

Similarly to~\cite{Zhang_15}, the gradient is computed at $t = 1$, $\nabla_{v_1} E(v_0) = \lambda(1) \cdot \nabla m(1)$
and integrated backward using the reduced adjoint Jacobi field equations (Equation~\ref{eq:AdjointJacobiField})
to obtain $\nabla_{v_0} E(v_0)$.
With this approach, the integration of the adjoint equation is not needed.

The second-order variations of the augmented Lagrangian on $w$ and $\lambda$ yield the incremental EPDiff and 
incremental state equations, needed for the computation of the Hessian-vector product. Thus,
\begin{eqnarray}
 \partial_t \delta v + ad_{\delta v}^\dagger v + ad_v^\dagger \delta v= 0 \textnormal{ in } \Omega \times (0,1] \\
 \partial_t \delta m + \nabla \delta m \cdot v + \nabla m \cdot \delta v = 0 \textnormal{ in } \Omega \times (0,1]
\end{eqnarray}
\noindent with initial conditions $\delta v(0) = 0$ and $\delta m(0) = 0$.

The Hessian-vector product $H_{v_0} E( \delta {v_0}^n)$ is computed from the Hessian-vector product at $t = 1$, which is integrated backward
using the reduced incremental adjoint Jacobi field equations
\begin{eqnarray}
 \partial_t \delta U + ad_{\delta v}^\dagger U + ad_v^\dagger \delta U= 0 \textnormal{ in } \Omega \times [0,1) \\
 \partial_t \delta w + \delta U - ad_{\delta v} w - ad_v \delta w + ad_{\delta w}^\dagger v + ad_{w}^\dagger \delta v = 0 \textnormal{ in } \Omega \times [0,1)
\end{eqnarray}

\noindent with initial conditions $\delta U(1) = \delta \lambda(1) \cdot \nabla m(1) + \lambda(1) \cdot \nabla \delta m(1)$, 
and $\delta w(1) = 0$. The proposed Gauss-Newton approximation drops $\lambda(1) \cdot \nabla \delta m(1)$ from the expression of $\delta U(1)$, 
and the $ad_v$ and $ad_w$ terms from the incremental adjoint Jacobi field equation on $\delta w$.

The minimization using a second-order inexact Gauss-Newton-Krylov method yields to the update equation
\begin{equation}
\label{eq:Newton}
 v_0^{n+1} = v_0^n - \epsilon \delta v_0^n,
\end{equation}
\noindent where $ \delta v_0^n$ is computed from PCG on the system 
\begin{equation}
\label{eq:KKT}
H_{v_0} E( \delta {v_0}^n) = - \nabla_{v_0} E(v_0^n).
\end{equation}
\noindent We also consider CG with the gradient and the Hessian computed on $V$ instead of $L^2$.

\section{Results}
\label{sec:Results} 

In this section, we evaluate the performance of our proposed method.
As a baseline for the evaluation, we include the results obtained by the methods most related to our work: 
Mang et al. method with the stationary and non-stationary parameterizations~\cite{Mang_15}, and 
Zhang et al. method in the spatial domain~\cite{Zhang_15}.

The experiments have been conducted on the Non-rigid Image Registration Evaluation Project database (NIREP).
In this work, images were resampled into volumes of size $180 \times 210 \times 180$. 
Registration was carried out from the first subject to every other subject in the database, yielding 15 registrations for each method.
The experiments were run on an NVidia GeForce GTX 1080 ti with 11 GBS of video memory and an Intel Core i7 with
64 GBS of DDR3 RAM.
We developed a hybrid Matlab 2017a CPU-GPU implementation of the methods using Cuda $8.0$.

\begin{figure*} [!t]
\centering
\scriptsize

\begin{tabular}{cc}
\includegraphics[width = 0.35 \textwidth]{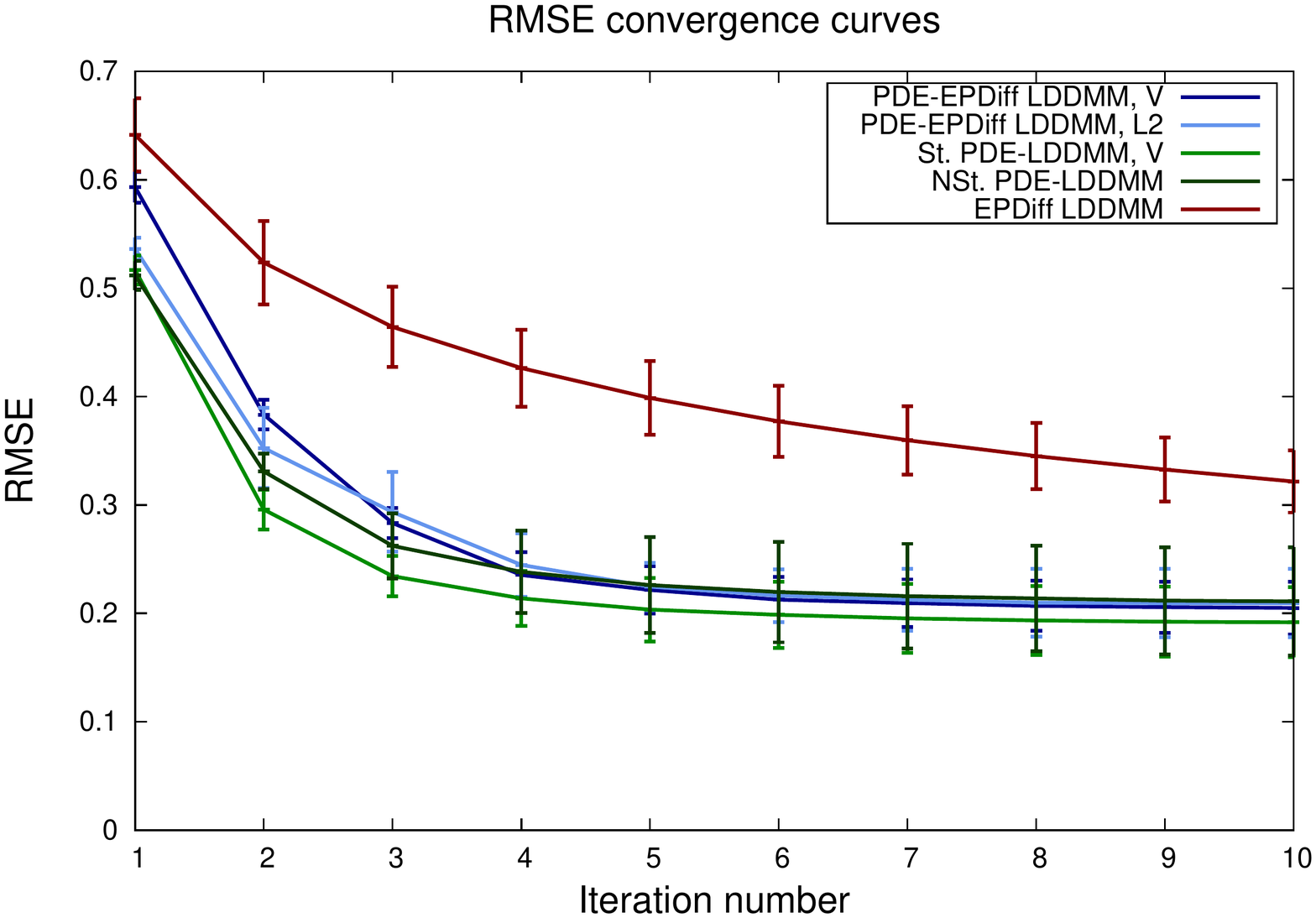} 
&
\includegraphics[angle = 0, width = 0.35 \textwidth]{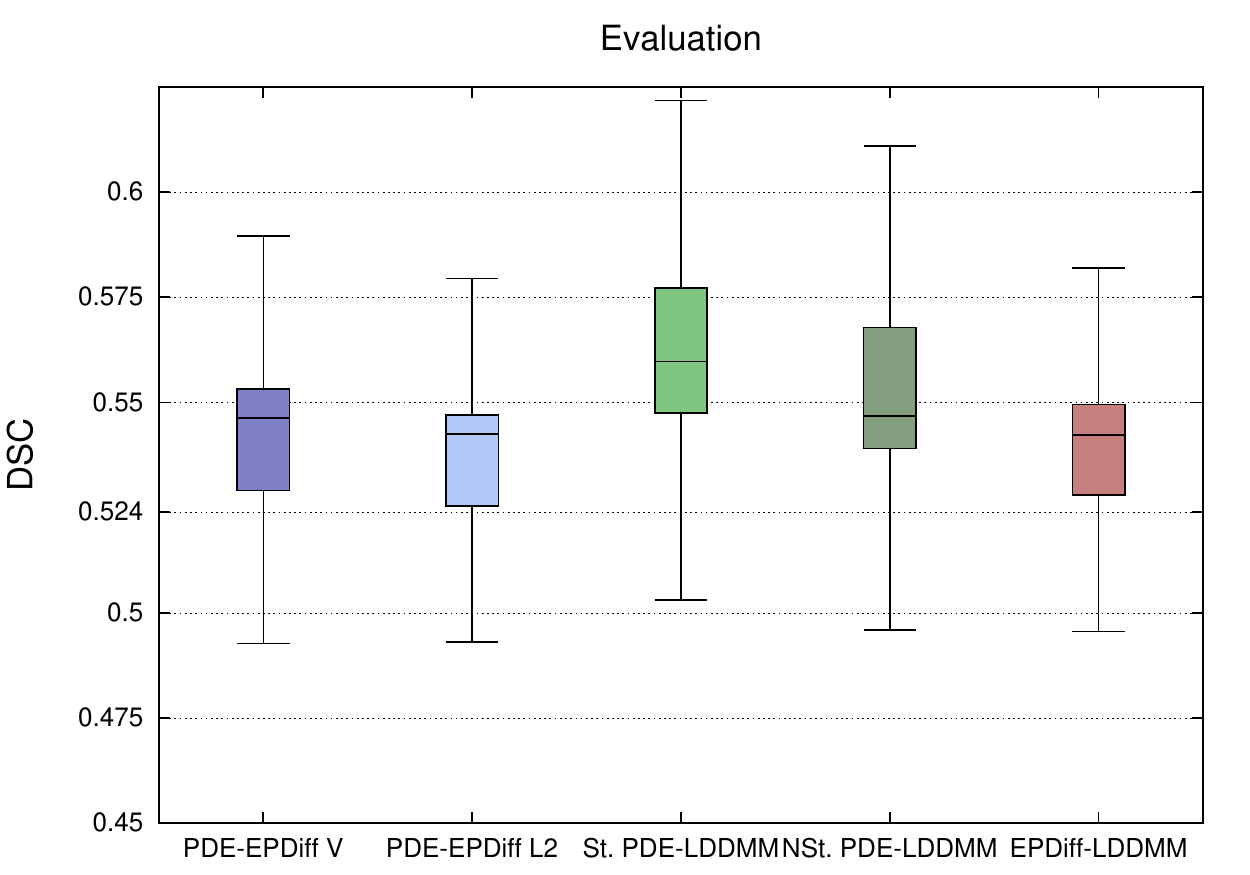} 
\end{tabular}
\caption{ \small Left, overall RMSE convergence curves. Right, mean and standard deviation of the DSC
coefficient obtained in the evaluation with the manual segmentations. } 
\label{fig:ConvergenceCurves}
\end{figure*}

Table~\ref{table:RMSJacs} shows, averaged by the number of experiments, 
the image similarity error $MSE_{rel} = \Vert m(1) - I_1\Vert_{L^2}^2 / \Vert I_0 - I_1 \Vert_{L^2}^2$, 
the relative gradient magnitude $\Vert g \Vert_{\infty, rel} = \frac{\Vert \nabla_v E(v^n) \Vert_\infty}{\Vert \nabla_v E(v^0) \Vert_\infty}$,
the extrema of the Jacobian determinant 
and the number of inner PCG iterations obtained after 10 iterations. 
Figure~\ref{fig:ConvergenceCurves} left shows the mean and standard deviation of the RMSE convergence curves 
obtained in the ten iterations.


\begin{table*}[!tpb]
\begin{center}
\scriptsize
\begin{tabular}{|c|c|c|c|c|c|c|c|}
\hline
Method & Opt. & $RMSE_{rel}$ & $\Vert g\Vert_{\infty,{rel}}$ & $\max(J(\phi_{1}^v)^{-1})$ & $\min(J(\phi_{1}^v)^{-1})$ & PCG iter \\
\hline
PDE-EPDiff LDDMM & GN, $V$ & 20.24 $\pm$ 2.05 & 0.12 $\pm$ 0.08 & 3.30 $\pm$ 1.12 & 0.11 $\pm$ 0.04 & 24.53 $\pm$ 4.14 \\
PDE-EPDiff LDDMM & GN, $L^2$ & 21.69 $\pm$ 2.87 & 0.28 $\pm$ 0.11 & 3.52 $\pm$ 1.14 & 0.11 $\pm$ 0.05 & 16.87 $\pm$ 3.50 \\
St. PDE-LDDMM & GN, $L^2$ & 19.18 $\pm$ 3.34 & 0.08 $\pm$ 0.05 & 3.71 $\pm$ 0.56 & 0.15 $\pm$ 0.04 & 41.00 $\pm$ 13.02 \\
NSt. PDE-LDDMM & GN, $L^2$ & 21.11 $\pm$ 5.18 & 0.20 $\pm$ 0.11 & 3.63 $\pm$ 0.83 & 0.14 $\pm$ 0.04 & 25.66 $\pm$ 14.86 \\
EPDiff LDDMM & GD, $V$   & 26.32 $\pm$ 2.23 & 0.07 $\pm$ 0.02 & 2.00 $\pm$ 0.17 & 0.28 $\pm$ 0.04 & 50.00 $\pm$ 0.00  \\
\hline
\end{tabular}
\\
\hspace{0.1 cm}
\caption{ 
Mean and standard deviation of the image similarity errors after registration, 
the relative gradient magnitude, 
the Jacobian determinant extrema associated with the transformation $(\phi_{1}^v)^{-1}$, 
and the number of PCG iterations.
In Opt. column, GN and GD stand for Gauss-Newton and gradient descent optimization, respectively.}
\label{table:RMSJacs}
\end{center}
\end{table*}




Figure~\ref{fig:ConvergenceCurves} right shows the accuracy of the registration results for template-based segmentation. 
We use the manual segmentations of the anatomical structures provided with the NIREP database as a 
gold standard and Dice Similarity Coefficient (DSC) is selected as the performance metric similarly to~\cite{Ou_14}. 

\begin{figure*} [!h]
\centering
\scriptsize
\begin{tabular}{|c|cc|cc|}
\hline
\includegraphics[angle = 0, width=2.0 cm, height = 2.0 cm]{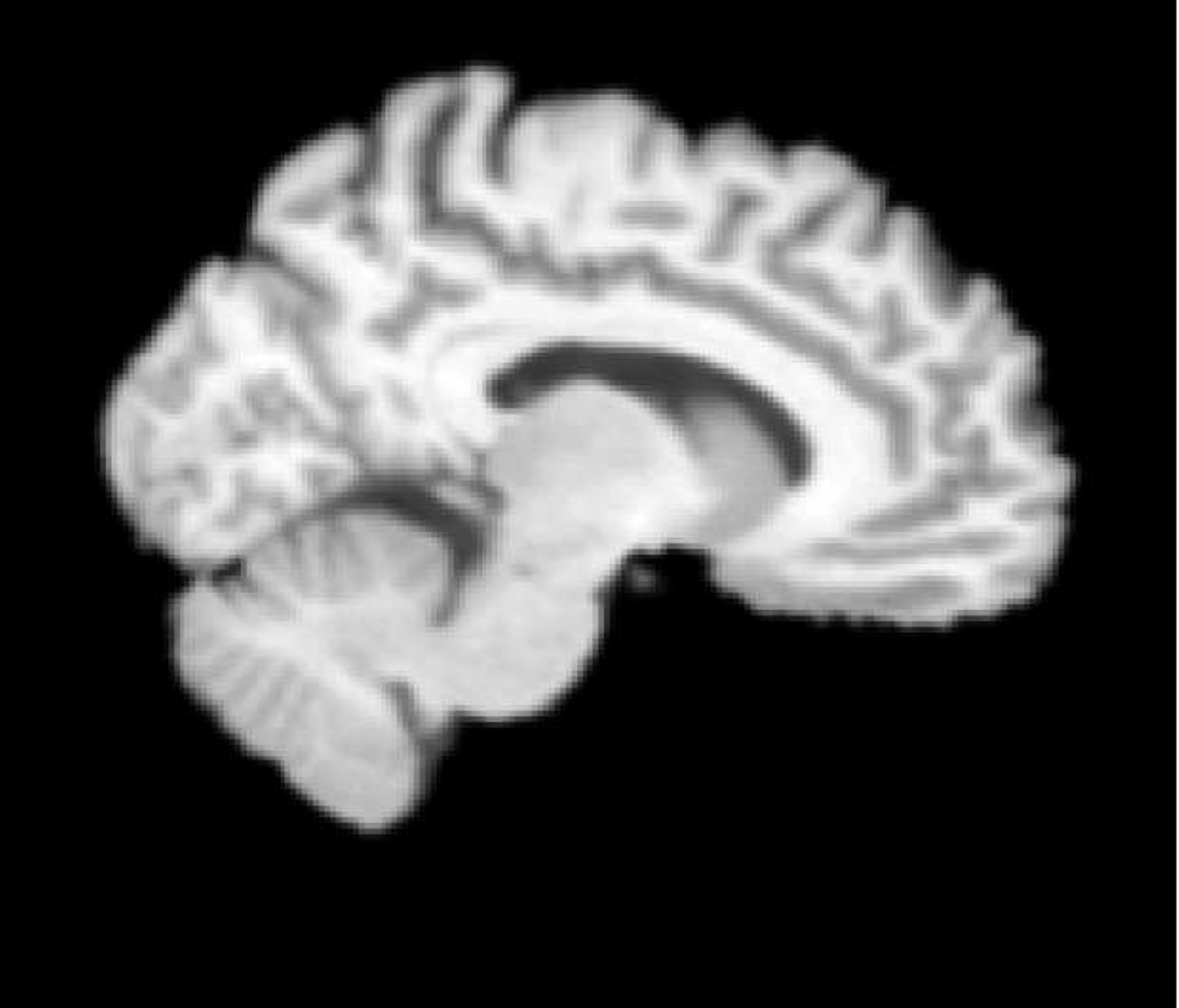} 
&
\includegraphics[angle = 0, width=2.0 cm, height = 2.0 cm]{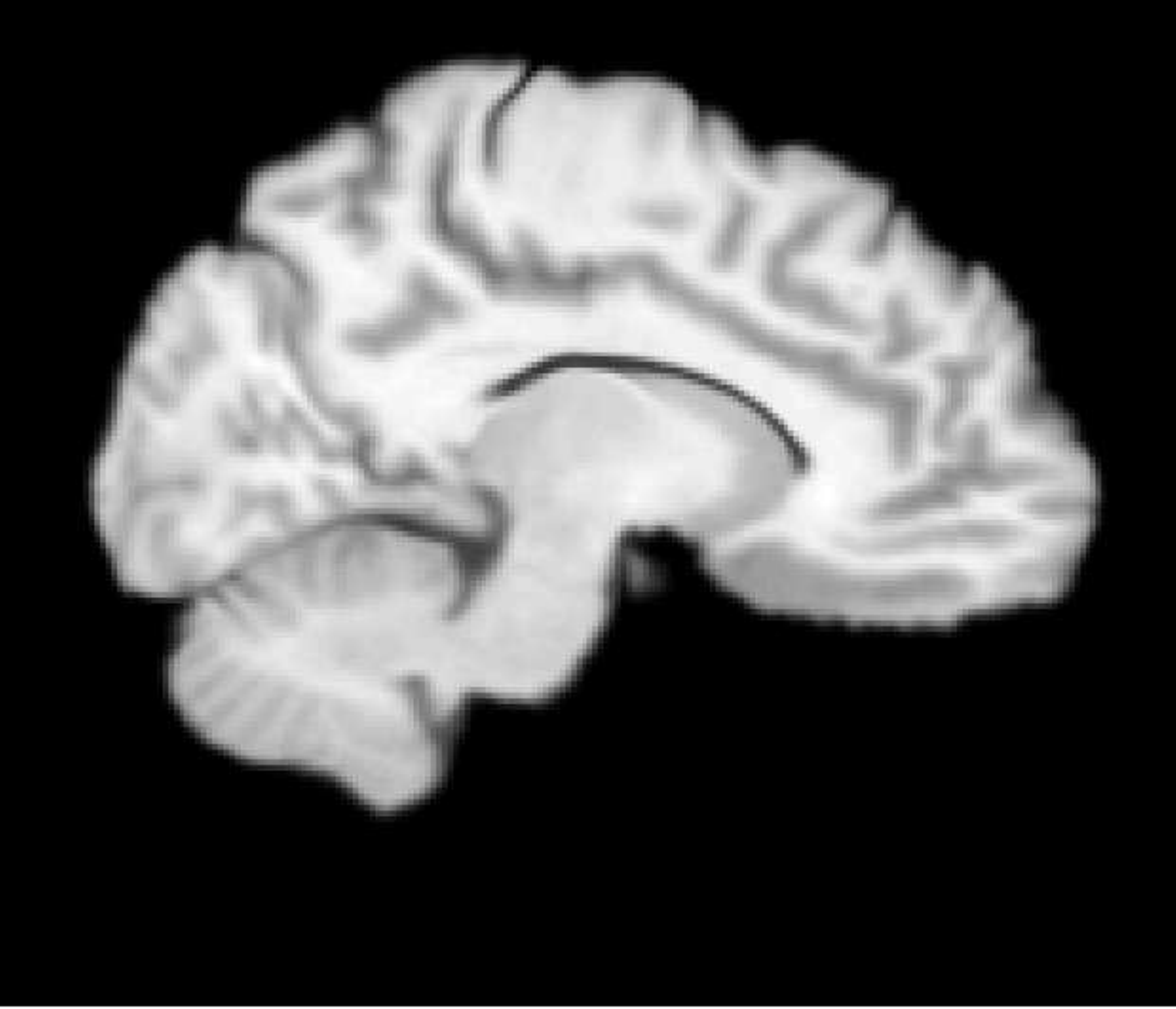} 
&
\includegraphics[angle = 0, width=2.0 cm, height = 2.0 cm]{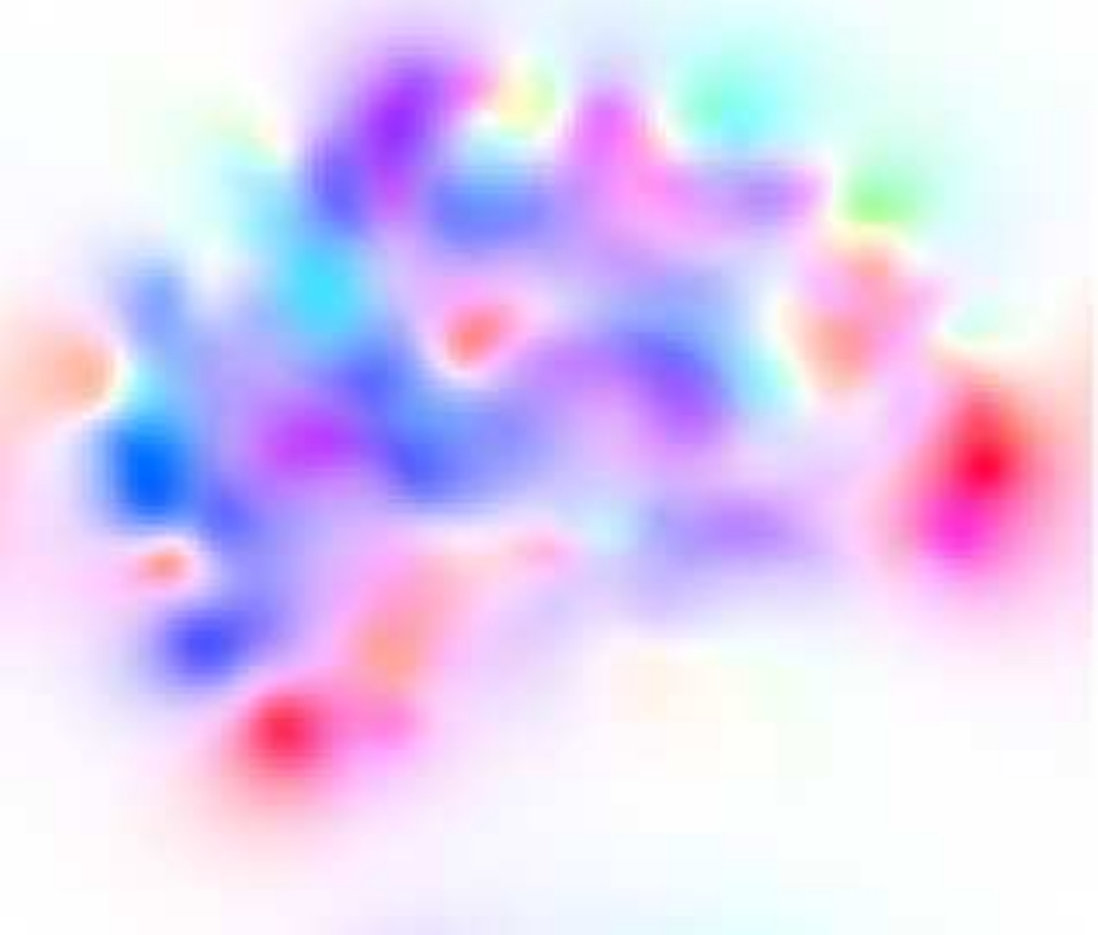} 
&
\includegraphics[angle = 0, width=2.0 cm, height = 2.0 cm]{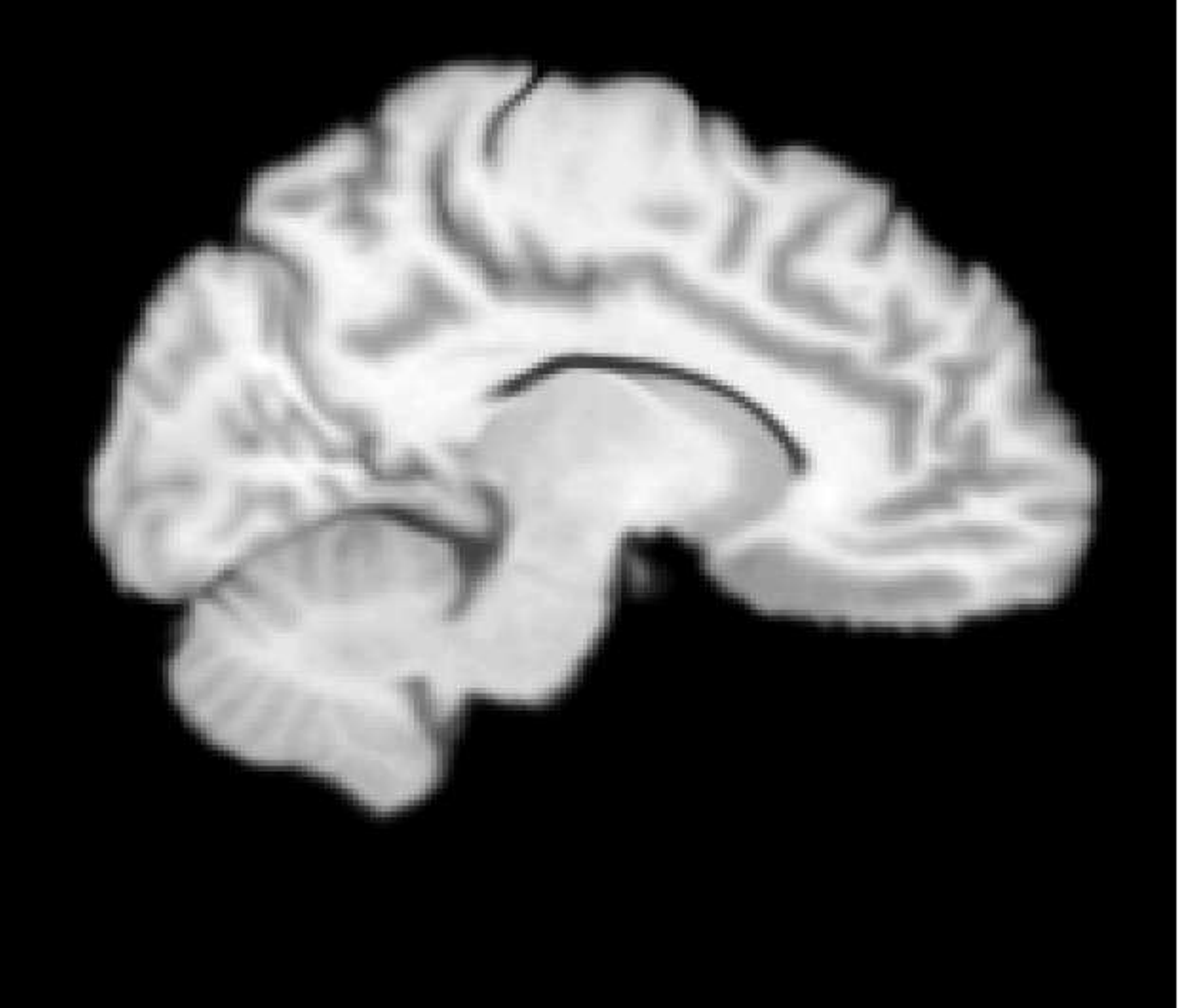} 
&
\includegraphics[angle = 0, width=2.0 cm, height = 2.0 cm]{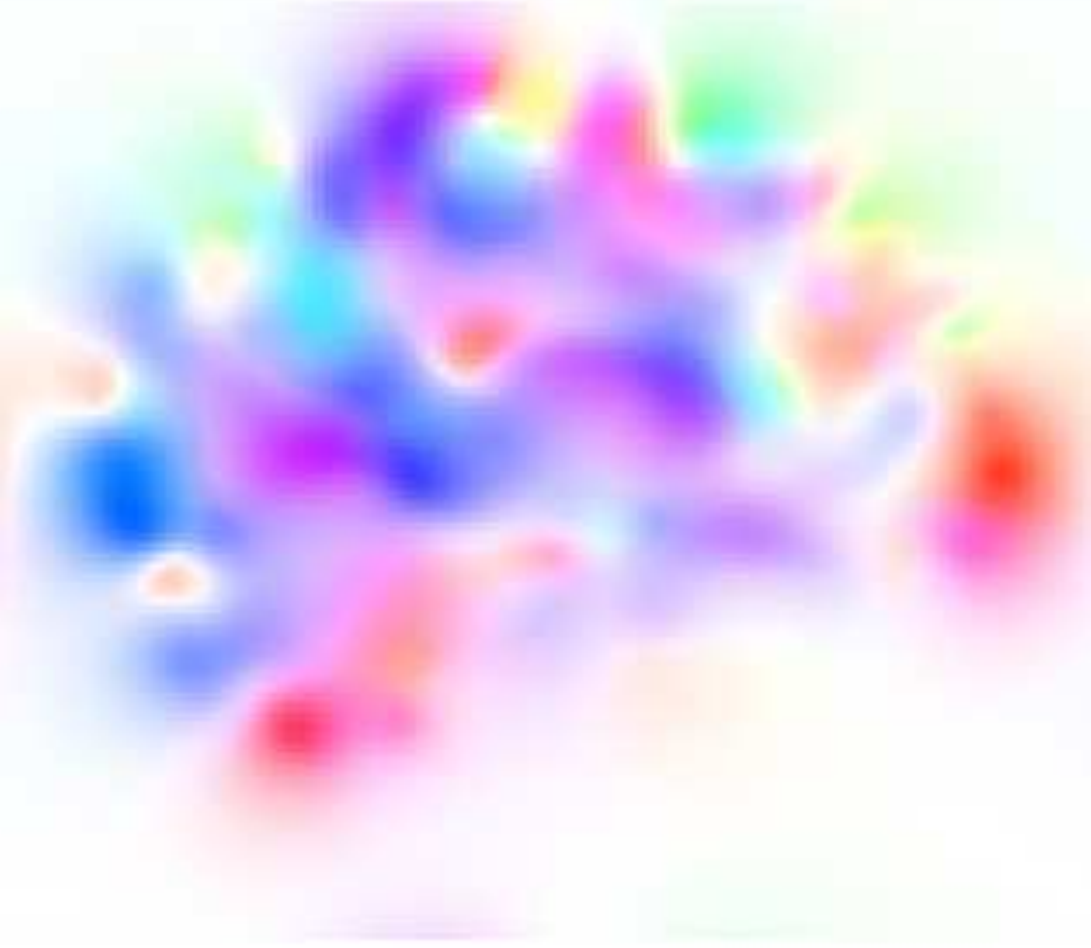} 
\\
source & PDE-EPDiff, V & &  St. PDE-LDDMM & 
\\
\includegraphics[angle = 0, width=2.0 cm, height = 2.0 cm]{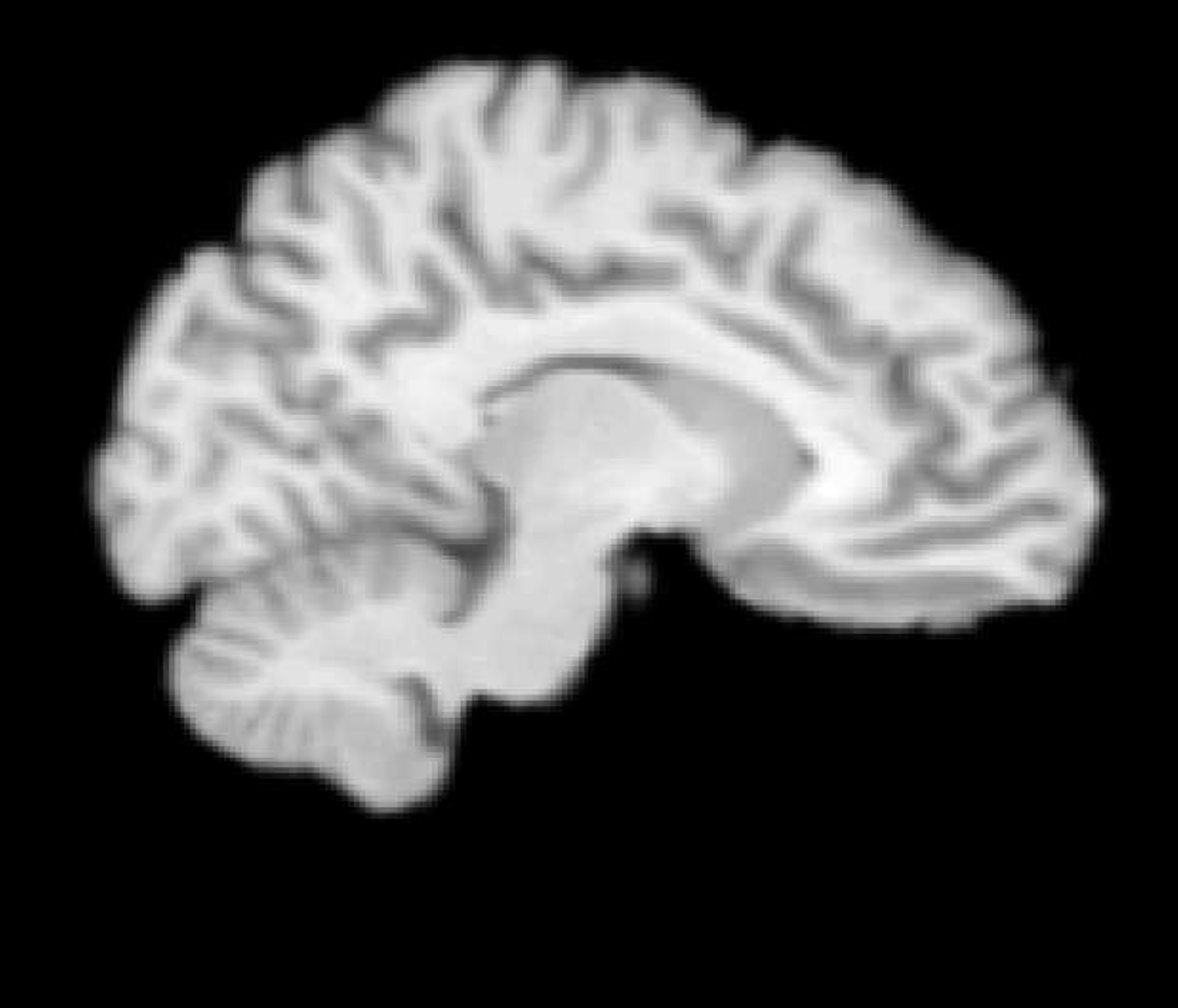} 
&
\includegraphics[angle = 0, width=2.0 cm, height = 2.0 cm]{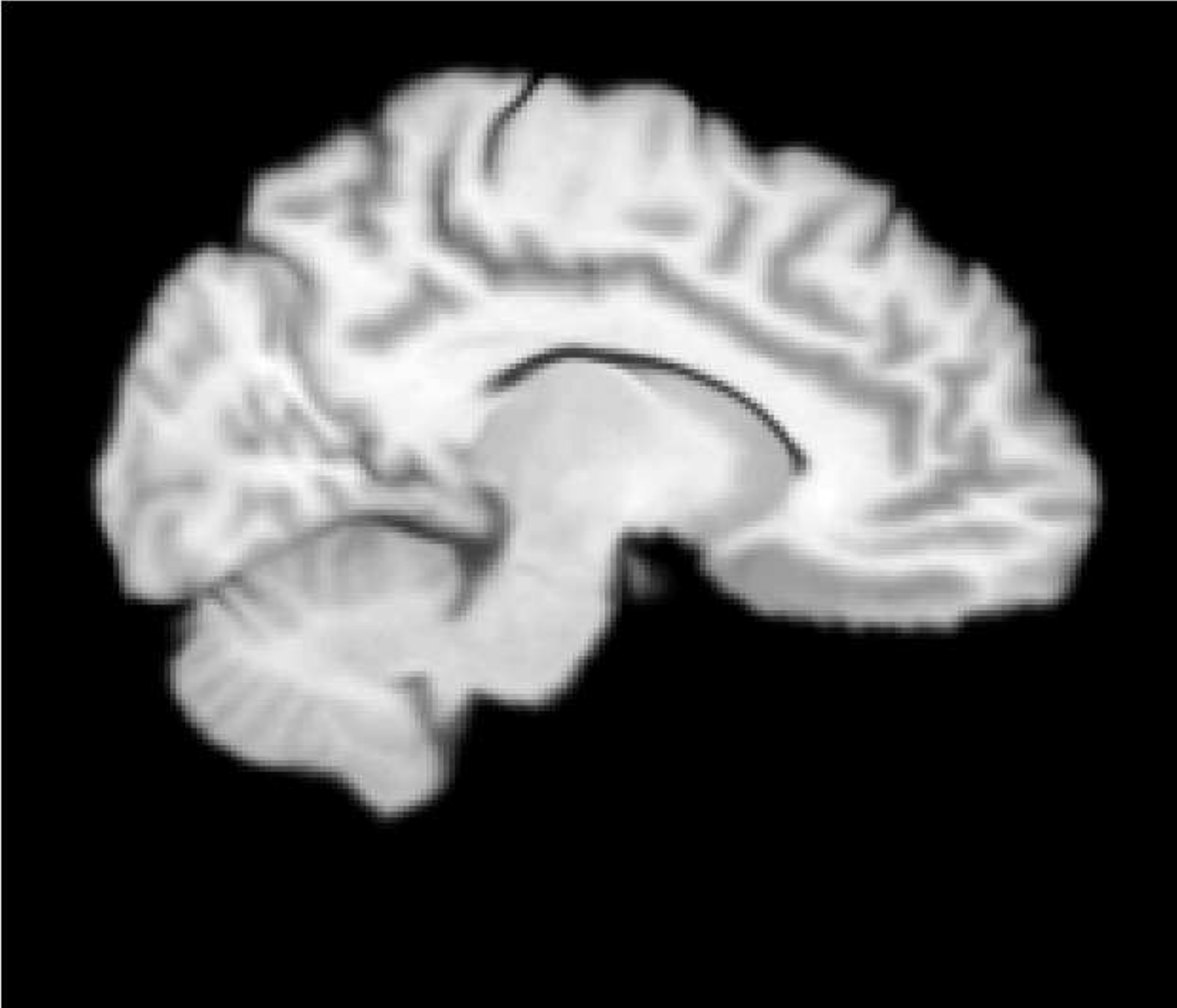} 
&
\includegraphics[angle = 0, width=2.0 cm, height = 2.0 cm]{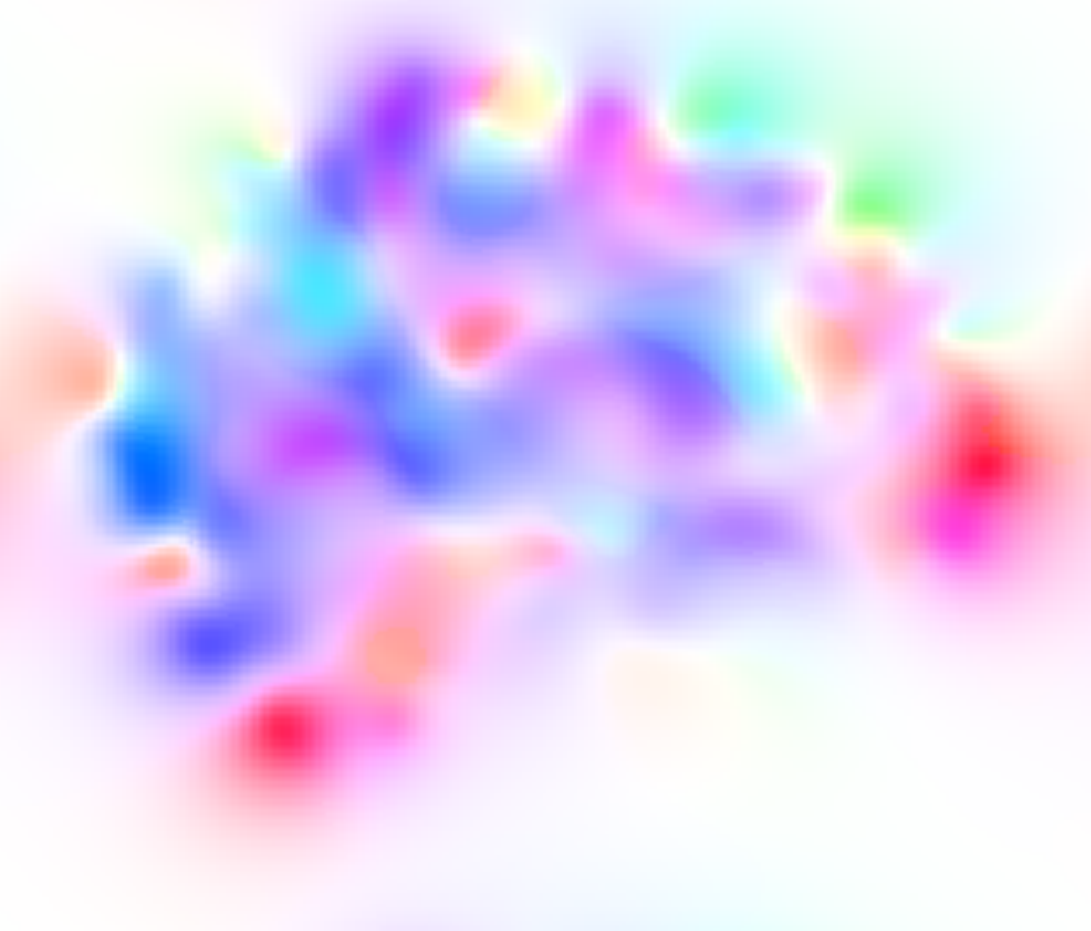} 
&
\includegraphics[angle = 0, width=2.0 cm, height = 2.0 cm]{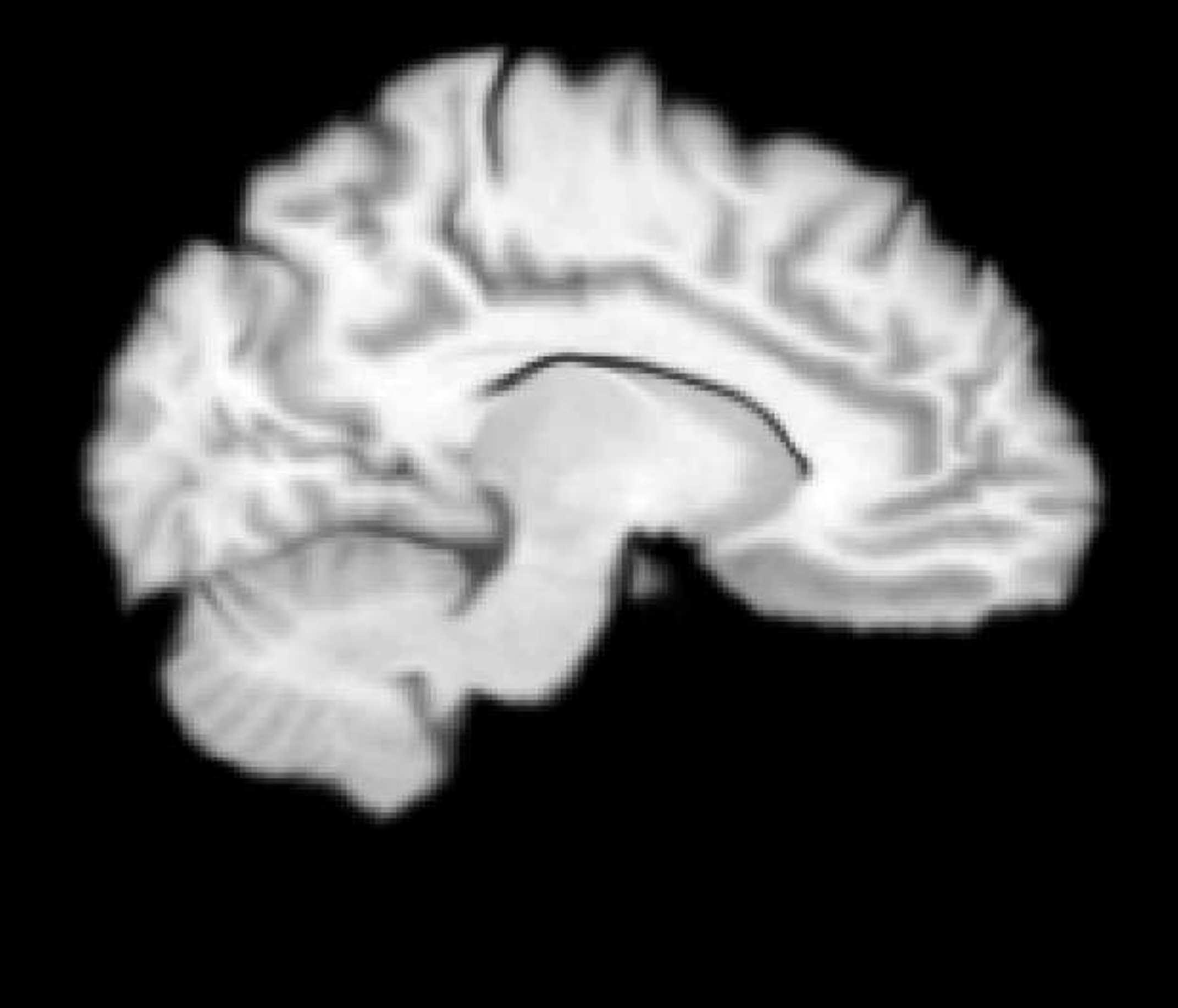} 
&
\includegraphics[angle = 0, width=2.0 cm, height = 2.0 cm]{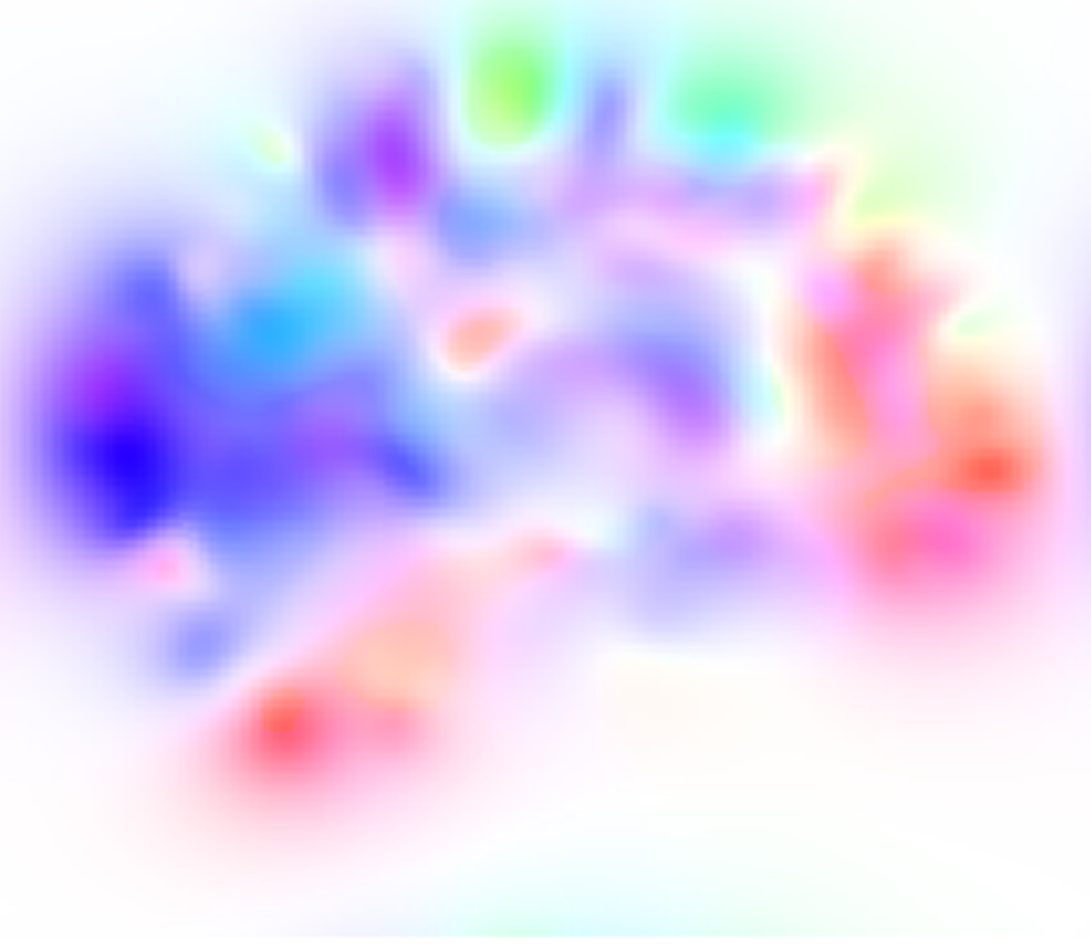} 
\\
target & PDE-EPDiff, $L^2$ & & EPDiff-LDDMM & \\
\hline
\end{tabular}
\caption{ \small Image registration results. Sagittal view of the warped sources and the 
initial velocity fields for the methods considered in the comparison.} 
\label{fig:3DResults}
\end{figure*}

Finally, for a qualitative assessment of the proposed registration method, we show the registration results 
in a selected experiment representative of a difficult deformable registration problem.
Figure~\ref{fig:3DResults} shows the warped images and the initial velocity fields obtained in the selected experiment.


\section{Discussion and Conclusions}
\label{sec:Conclusions}

In this work, we have proposed a method for PDE-constrained LDDMM parameterized in the space of initial velocity 
fields. 
Our method can be regarded as an extension of the method in~\cite{Mang_16} that allows obtaining geodesics.
We propose to perform the derivations of the gradient and Hessian-vector products on the final velocity field 
and transport the computations backward using the adjoint and the incremental adjoint Jacobi equations.
This way the complex dependence on the initial velocity field is skipped in the derivations.
We also avoid the computation of the adjoint equation and its incremental counterpart that has been recently 
identified as a subtle problem~\cite{Mang_17}.

The method has been compared and evaluated with respect to benchmark PDE-constrained LDDMM and EPDiff-LDDMM methods.
Gauss-Newton-Krylov optimization provided a higher rate of convergence than gradient-descent as expected.
The evaluation in the NIREP database showed competitive performance with respect to the benchmark methods
despite the more restrictive model imposed by the EPDiff equation on the velocity fields.

The major drawback of the proposed method is the large memory load inherent to PDE-constrained LDDMM methods.
We will deal with this problem in future work.


%
%
%

\section*{Acknowledgements}
The author would like to give special thanks to Wen Mei Hwu from the University of Illinois 
for interesting ideas in the GPU implementation of the methods.
This work was partially supported by Spanish research grant TIN2016-80347-R.

\bibliographystyle{splncs}
\bibliography{abbsmall.bib,Diffeo.bib,OpticalFlow.bib}

\end{document}